\title{Non-backtracking random walks mix faster}
\author{{Noga Alon\thanks{School of Mathematics, Raymond and Beverly Sackler
Faculty of Exact Sciences, Tel Aviv University, Tel Aviv, 69978,
Israel. Email: nogaa@tau.ac.il. Research supported in part by the
Israel Science Foundation, by a USA-Israeli BSF grant, and by the
Hermann Minkowski Minerva Center for Geometry at Tel Aviv University.}}
\quad {Itai Benjamini \thanks{Weizmann Institute, Rehovot, 76100,
Israel. Email: itai.benjamini@weizmann.ac.il}} \quad {Eyal Lubetzky
\thanks{ School of Computer Science, Raymond and Beverly
Sackler Faculty of Exact Sciences, Tel Aviv University, Tel Aviv,
69978, Israel. Email: lubetzky@tau.ac.il. Research partially
supported by a Charles Clore Foundation Fellowship.}} \quad {Sasha
Sodin
\thanks{ School of Mathematics, Raymond and Beverly Sackler Faculty
of Exact Sciences, Tel Aviv University, Tel Aviv, 69978, Israel.
Email: sodinale@tau.ac.il.}}}

\documentclass [11pt]{article}
\usepackage[]{amsmath,amssymb,amsfonts,latexsym,amsthm,enumerate}
\usepackage[]{graphicx}

\newtheorem{theorem}{Theorem}[section]
\newtheorem{lemma}[theorem]{Lemma}
\newtheorem{claim}[theorem]{Claim}
\newtheorem*{definition}{Definition}

\newtheorem{corollary}[theorem]{Corollary}

\renewcommand{\epsilon}{\varepsilon}

\newcommand{\cov}{\operatorname{Cov}}
\newcommand{\sign}{\operatorname{sign}}
\newcommand{\var}{\operatorname{Var}}

\newtheoremstyle{upright}%
        {8pt plus2pt minus4pt}%
        {8pt plus2pt minus4pt}%
        {\upshape}%
        {}%
        {\bfseries}%
        {:}%
        {1em}%
        {}%
\theoremstyle{upright}
\newtheorem{remark}[theorem]{Remark}

\newcommand{\ignore}[1]{}

\begin{document}
\maketitle

\begin{abstract}We compute the mixing rate of a non-backtracking random walk on a regular
expander. Using some properties of Chebyshev
polynomials of the second kind, we show that this rate may be up
to twice as fast as the mixing rate of the simple random walk. The closer the
expander is to a Ramanujan graph, the higher the ratio between the
above two mixing rates is.

As an application, we show that if $G$ is a high-girth regular
expander on $n$ vertices, then a typical non-backtracking random
walk of length $n$ on $G$ does not visit a vertex more than
$(1+o(1))\frac{\log n}{\log\log n}$ times, and this result is tight.
In this sense, the multi-set of visited vertices is analogous to the
result of throwing $n$ balls to $n$ bins uniformly, in contrast to
the simple random walk on $G$, which almost surely visits some
vertex $\Omega(\log n)$ times.
\end{abstract}

\section{Introduction}
\subsection{Background and definitions}
 Let $G=(V,E)$ be an undirected graph. A random walk of length $k$
on $G$, from some given vertex $w_0 \in V$, is a uniformly chosen
member of:
$$\mathcal{W}^{(k)} = \{
(w_0,w_1,\ldots,w_k) ~:~w_t \in V,~w_{t-1} w_t \in E~\mbox{ for all
} t\in[k]\}~.$$ Equivalently, such a walk is a finite Markov chain
$\mathcal{M}=(X_0,\ldots,X_k)$ on the state space $V$, where
$X_0=w_0$ and the transition probabilities are $P_{u v} = \Pr[ X_i =
v ~|~ X_{i-1} = u] = \mathbf{1}_{\{u v \in E\}} /\deg(u)$. For
further information on Markov chains, see, e.g, \cite{Keilson},
\cite{Norris}.

The extensive study of random walks on graphs was motivated by the
following useful property, which we first state informally. While
the random walk is simple to analyze and to implement in many
frameworks, it ``mixes'' in $G$ after a relatively small number of
steps, provided $G$ satisfies some natural requirements. Thus,
the random walk provides an efficient method of sampling the graph
vertices, a fact which has many applications in Theoretical and
Applied Computer Science. See \cite{Lovasz} for a survey on the
subject.

The following facts are well known (see, e.g., \cite{Lovasz},
\cite{LovaszWinkler}, \cite{Sinclair}). If $G$ is a connected and
non-bipartite undirected graph, then the Markov chain $\mathcal{M}$,
corresponding to the random walk on $G$, is irreducible and
aperiodic. In this case, $\mathcal{M}$ converges to a unique
stationary distribution, $\pi$, regardless of its starting position,
where $\pi(u) = \frac{\deg(u)}{2|E|}$. The {\em mixing rate} of the
random walk on $G$ measures how fast $\mathcal{M}$ converges to the
stationary distribution, and is defined as follows:
\begin{equation}
  \label{eq-mixing-rate-def}
  \rho = \rho(G) = \limsup_{k\to\infty} \max_{u,v\in V} \left| P_{u v}^{(k)} -
  \pi(v) \right| ^{1/k}~,
\end{equation}
where $P_{u v}^{(k)} = \Pr[X_{t+k}= v~|~X_t = u]$. The notion of
{\em mixing time}, the number of steps it takes $\mathcal{M}$ to get
``sufficiently close'' to $\pi$, has several commonly used
definitions, and for each of these definitions there are lower and
upper bounds as a function of $\rho$ and $n$. For instance, letting
$P^{(k)}_u$ denote the distribution of $\mathcal{M}$ at time $k$
given that $X_0=u$, one may define the mixing time $\tau_\epsilon$
as the minimal number of steps it takes $P^{(k)}_u$ and $\pi$ to be
at most $\epsilon$-far in terms of their total variation distance,
maximized over all vertices $u \in V$.

An important special case of the above is the one where the graph
$G$ is regular. In this case, the stationary distribution $\pi$ is
the uniform distribution, being an eigenvector of the transition
probabilities matrix $P=A / d$ (where  $A$ is the adjacency matrix
of the graph and $d$ is its regularity degree).  Hence, whenever $G$
is connected and non-bipartite, the random walk eventually approximates the
uniform distribution. As we next
specify, these sufficient and necessary conditions, required for the
random walk on $G$ to mix, are determined by the spectrum of $G$.

 Let $G$ be a $d$-regular graph. The eigenvalues of $G$, that
is, the eigenvalues of its (symmetric) adjacency matrix are
$d=\lambda_1 \geq \lambda_2 \geq\ldots \geq \lambda_n$, and
$|\lambda_i|\leq d$ for all $i$ (by the Perron-Frobenius Theorem).
The multiplicity of the eigenvalue $d$ is equal to the number of
connected components of $G$, and $\lambda_n = -d$ iff $G$ is
bipartite (proofs of these well known facts can be found, for
instance, in \cite{AlgabraicGraphTheory}). Therefore, whenever $G$
is $d$-regular, the conditions that $G$ should be connected and
non-bipartite become equivalent to requiring that
$\lambda=\max\{\lambda_2,|\lambda_n|\}$ would satisfy $\lambda<d$.
Define the following:
\begin{definition}
  An $(n,d,\lambda)$-graph, for some integer $d$ and some
$\lambda < d$, is a $d$-regular graph on $n$ vertices whose second
largest eigenvalue in absolute value is $\lambda$.
\end{definition}
This notion was introduced by the first author in the 80's,
motivated by the fact that if $\lambda$ is much smaller than $d$, then
the graph has strong pseudo-random properties. We mention a few of the
properties of these graphs, and refer the readers to \cite{KS} for
an extensive survey of the subject. Let $G=(V,E)$ denote an
$(n,d,\lambda)$-graph. First, the behavior of $G$ resembles that of
a random graph of edge density $d/n$ in the following sense: if $A,B
$ are (not necessarily disjoint) subsets of vertices, then $\left|
e(A,B) - \frac{d}{n} |A||B| \right| \leq \lambda\sqrt{ |A| |B|}~,$
where $e(A,B)$ denotes the number of ordered pairs $\{(a,b):a \in A,
b\in B, a b \in E\}$ (see \cite{ProbMethod}, Corollary 9.25). In
other words, every two sets of vertices $A,B$ have roughly the ``right''
number of edges between them. Second, the expansion property of $G$
is closely related to the eigenvalue gap $d-\lambda$, as stated
next. Defining the vertex boundary of $X$, $\delta X$, as the set of
neighbors of $X$ in $V \setminus X$, it is known that $|\delta X|
\geq \frac{2(d-\lambda)}{3d-2\lambda}|X|$ for all sets $X$ of size
at most $n/2$ (\cite{AlonMilman}). Conversely, if
\begin{equation}\label{eq-expander-def}
|\delta X| \geq c |X|~\mbox{ for some $c > 0$ and all }X\subset
V~,|X|\leq \frac{n}{2}~,\end{equation} then $d-\lambda \geq
c^2/(4+2c^2)$, implying a discrete version of Cheeger's inequality
(\cite{Alon}). A graph satisfying \eqref{eq-expander-def} with $c$ bounded
away from $0$ is commonly referred to as an {\em expander}, and according to this
definition $(n,d,\lambda)$-graphs with $\lambda$ bounded away from $d$
and regular expanders are very close notions.

In many applications of random walks on expanders, there is not much
sense in allowing the walk to backtrack, besides making the model
easier to understand and to analyze. A non-backtracking random walk
on an undirected graph $G$, is a walk which does not traverse the
same edge twice in a row.
In the first part of this paper, we determine the mixing rate of
non-backtracking random walks on expanders, using some properties of
Chebyshev polynomials of the second kind (the connection between
these polynomials and non-backtracking walks follows ideas
from \cite{McKay,LiSole}). We obtain that for $3 \leq d \leq n^{o(1)}$, the
mixing-rate of a non-backtracking random walk on an
$(n,d,\lambda)$-graph is at most the mixing-rate of a simple random
walk on the same graph. In fact, the ratio between the two may reach up to
$\frac{2(d-1)}{d}$, as formulated in the next Subsection.

Let $G$ be a $d$-regular expander. The following
definition of the mixing-time of a random walk on $G$ corresponds to
an $L_\infty$ distance of $\frac{1}{2n}$, as well as to a relative
pointwise distance (r.p.d.) of $\frac{1}{2}$, between $\pi$ and
$P^{(k)}_u$, for all $u\in V$:
\begin{equation}\label{eq-mixing-time-def}
\tau = \tau(G) = \min_t \left\{ \Big|P^{(k)}_{u v} - \frac{1}{n}
\Big| \leq \frac{1}{2n} ~\mbox{ for all $u,v\in V$ and $k \geq
t$}\right\}~.
\end{equation}
As $G$ is a regular expander, $\tau = \Theta(\log
n)$. Notice that sampling the position of the random walk
$\mathcal{M}$ at time-points, which are at least $\tau$-apart, gives
a more or less independent and uniformly distributed set of
vertices. On the other hand, a set of vertices sampled at
constant-distance time-points is clearly very much dependent. As we
next show, there is a special interest in the distribution of the
set of vertices along $\Theta(n)$ consecutive steps of the random
walk.

An example of this is the amplification of randomized algorithms
(such as the Rabin-Miller primality testing algorithm).
Let $\mathcal{A}$ denote such an algorithm which uses $\log n$ random
bits; the naive parallel repetition of $\mathcal{A}$ spends
$\Theta(n \log n)$ bits in order to reduce the error probability
exponentially in $n$. It is well known that the probability that a
random walk of length $k$ avoids a given set of vertices of constant
proportion, decreases exponentially with $k$ (see, e.g., \cite{ProbMethod},
Corollary 9.28). Therefore, if $G$ is a regular expander of fixed
degree, feeding the positions of a random walk of length $\Theta(n)$
on $G$ as the random seeds for the algorithm, reduces the error
probability of the algorithm exponentially, using only $\Theta(n)$
random bits.

In the above application of conserving randomness when amplifying
randomized algorithms, our concern was the probability that a random walk
of length $n$ misses a large given set of vertices. Instead, in load
balancing applications, the concern is the maximal number of times
that a random walk of length $n$ visits a vertex. This corresponds
to the classical {\em balls and bins} paradigm (see, e.g,
\cite{Feller}, \cite{JohnsonKotz}), which discusses the result of
throwing $n$ balls to $n$ bins, independently and uniformly at
random. In the balls and bins experiment, the bin with the largest
number of balls typically contains $(1+o(1))\frac{\log n}{\log\log
n}$ balls (see \cite{Gonnet}).

As we later show, the random walk is unsuitable for conserving
randomness in this case, as a typical random walk of length $n$ has
a maximal load of $\Omega(\log n)$. As an application for
non-backtracking random walks, we show here that the maximal number
of times that such a walk of length $n$ visits a vertex, is
$(1+o(1))\log n / \log\log n$ times on high girth expanders
with $n$ vertices.

Throughout the paper, we say that an event, which is defined for an
infinite series of graphs, occurs {\em with high probability}, or
{\em almost surely}, or that {\em almost every} graph of an infinite
series of graphs satisfies some property, if the probability for the
corresponding event tends to $1$ as the number of vertices tends to
infinity. Unless stated otherwise, all logarithms are in the natural
basis.

\subsection{Main results} Let $G=(V,E)$ denote
an undirected graph. Define a {\em non-backtracking random walk} of
length $k$ on $G$, from some given vertex $w_0 \in V$, as a
uniformly chosen member of:
$$\widetilde{\mathcal{W}}^{(k)} = \left\{
(w_0,w_1,\ldots,w_k) ~:~\begin{array}{c} w_t \in V~,~w_{t-1} w_t \in
E~\mbox{ for all } t\in[k]~,\\ w_{t-1}\neq w_{t+1}~\mbox{ for all
}t\in[k-1]\end{array}\right\}~.$$ Equivalently, a non-backtracking
random walk on $G$ from $w_0$ is a finite Markov chain
$\widetilde{\mathcal{M}}$, whose state space is
$\overrightarrow{E}$, the set of directed edges of $G$, taking each
edge in both orientations. The distribution of the initial state is
given by $\Pr[X_0 = (w_0,u)] = \mathbf{1}_{\{w_0 u \in
E\}}/\deg(w_0)$ (and $0$ elsewhere), and the transition
probabilities are $P_{(u,v),(v,w)} = \mathbf{1}_{\{u \neq w\}}
/(\deg(v)-1)$ (and $0$ elsewhere). If $G$ is $d$-regular, then the
transition probabilities matrix is double-stochastic, hence the
uniform distribution is a stationary distribution of
$\widetilde{\mathcal{M}}$. Notice that if $G$ is $2$-regular, then
it is a disjoint union of cycles, hence a non-backtracking random
walk on $G$ is periodic and does not converge to a stationary
distribution. We therefore require that $d \geq 3$, in addition to
the requirements that $G$ should be connected and non-bipartite, and
these necessary conditions prove to be sufficient for
$\widetilde{\mathcal{M}}$ to converge to the uniform distribution.

 Let $G=(V,E)$ denote an $(n,d,\lambda)$-graph for $d \geq 3$. Recalling \eqref{eq-mixing-rate-def},
define the mixing rate of a non-backtracking random walk on $G$ as:
\begin{equation}
  \label{eq-widetilde-rho-def}
  \widetilde{\rho}(G)=\limsup_{k\to\infty} \max_{u,v\in V} \Big|
\widetilde{P}_{u v}^{(k)} - \frac{1}{n} \Big| ^{1/k}~,
\end{equation} where
$\widetilde{P}_{u v}^{(k)}$ is the probability that a
non-backtracking random walk of length $k$ on $G$, which starts in
$u$, ends in $v$. The following theorem, proved in Section
\ref{sec::mixing}, determines the value of $\widetilde{\rho}$ in
this case:
\begin{theorem}\label{thm-1} Let $d \geq 3$ denote some
integer, and let $G$ be an $(n,d,\lambda)$-graph for some $\lambda <
d$. Define $\psi:[0,\infty)\to\mathbb{R}$ by:
\begin{equation}\label{eq-psi-def} \psi(x) = \left\{\begin{array}
  {ll}x+\sqrt{x^2-1} & \mbox{If }~x \geq 1~,\\
1 & \mbox{If }~0 \leq x \leq 1~.
\end{array}\right.\end{equation}
Then a non-backtracking random walk on $G$ converges to the uniform
distribution, and its mixing rate, $\widetilde{\rho}$, satisfies:
\begin{equation}\label{eq-rho-values}
  \widetilde{\rho}=
\psi\left(\frac{\lambda}{2\sqrt{d-1}}\right)/\sqrt{d-1}~.
\end{equation}
\end{theorem}

\begin{figure}
\centering
\begin{tabular}{cc}
\includegraphics[width=3in]{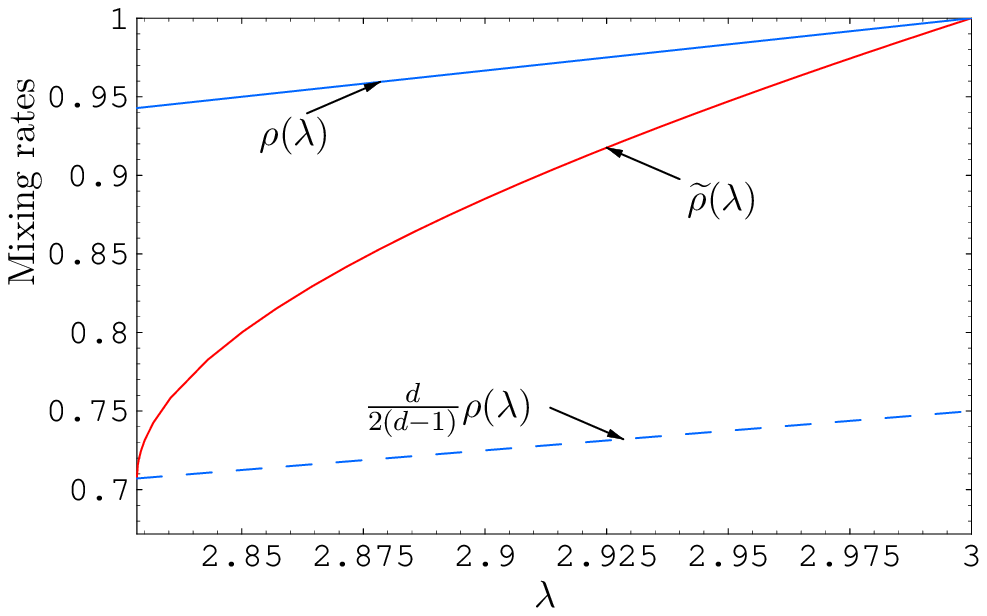} &
\includegraphics[width=3in]{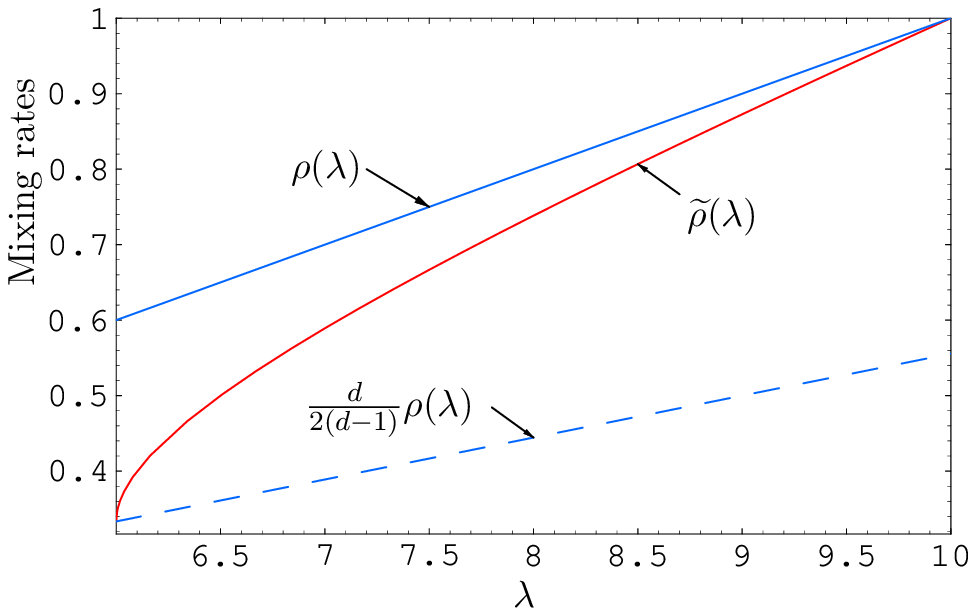} \\
\mbox{\small (a) $3$-regular graphs} & \mbox{\small (b) $10$-regular
graphs}
\end{tabular}
\caption{ Mixing rates of simple and non-backtracking random walks
on regular graphs. } \label{fig:mixing}
\end{figure}
It is well known (see, e.g., \cite{Lovasz}), that if $G$ is an
$(n,d,\lambda)$-graph, then the mixing-rate of the simple random
walk on $G$ is $\rho = \lambda / d$. As we state in Section
\ref{sec::mixing}, combining this with the properties of the
function $\psi$, defined in \eqref{eq-psi-def}, gives the inequality
$\widetilde{\rho} \leq \rho$, provided $d \leq n^{o(1)}$. The closer $\lambda$ is to
$2\sqrt{d-1}$ (that is, the closer the graph is to being Ramanujan),
the closer the ratio $\widetilde{\rho}/\rho$ is to
$\frac{d}{2(d-1)}$, as demonstrated in Figure \ref{fig:mixing}. This
is formulated in the following corollary:
\begin{corollary}\label{cor-mixing-ratio}
Let $G$ be a non-bipartite and connected $d$-regular graph on $n$
vertices, for some $d\geq 3$, and let $\rho$ and
$\widetilde{\rho}$ denote the mixing rates of simple and
non-backtracking random walks on $G$, respectively. The following
holds: let $\lambda$
be the second largest eigenvalue of $G$ in absolute value. If $\lambda \geq 2\sqrt{d-1}$, then
\begin{equation}\label{eq-mixing-rates-ratio} \frac{d}{2(d-1)}
\leq \frac{\widetilde{\rho}}{\rho} \leq 1~.\end{equation}
If $\lambda < 2\sqrt{d-1}$ and $d=n^{o(1)}$, then
$\widetilde{\rho} /\rho = \frac{d}{2(d-1)}+o(1)$,
 where the
$o(1)$-term tends to $0$ as $n\to\infty$.
\end{corollary}

In Section \ref{sec::balls-bins}, we discuss the maximal load of a
set of vertices along $n$ consecutive positions of a
non-backtracking random walk. The next theorem states that the
maximal number of times that such a walk on a regular expander of
high girth visits a vertex is equal to $(1+o(1))\frac{\log
n}{\log\log n}$, precisely the maximal load in the balls and bins
experiment.

\begin{theorem}\label{thm-2} Let $G$ be an $(n,d,\lambda)$ graph for some
fixed $d \geq 3$ and some fixed $\lambda < d$, whose girth is $g
\geq 10 \log_{d-1}\log n$. With high probability, the maximal number
of times that a non-backtracking random walk of length $n$ on $G$
visits a vertex is equal to $(1+o(1))\frac{\log n}{\log\log n}$.
\end{theorem}
Furthermore, the above requirement on the girth is essentially
tight: in Section \ref{sec::balls-bins} we show that, for all
$g=g(n)$, there are graphs as described in Theorem \ref{thm-2} with
girth $g$, for which the above maximal number of visits is
$\Omega(\frac{\log n}{g})$ almost surely.

The final section, Section \ref{sec::conclusion}, is devoted to
several open problems, further related to random walks on expanders
and to similar notions of conserving randomness.

\section{The mixing rate of a non-backtracking random
walk}\label{sec::mixing}
\noindent{\textbf{Proof of Theorem \ref{thm-1}.}
We begin with some preliminaries on Chebyshev polynomials; for
further information, see, e.g., \cite{Szego}. The Chebyshev
polynomials of the second kind, of degree $k\geq 0$, are
the following polynomials:
\begin{align} U_k(\cos\theta) &=
\frac{\sin\left((k+1)\theta\right)}{\sin \theta} \label{eq-cheb-def}
~.
\end{align}
Also, it is convenient to define $U_{-1}(x) \equiv 0$. The Chebyshev polynomials satisfy the following three-term recurrence
relation:
\begin{equation}
  \label{eq-cheb-recursion}
  U_{k+1}(x) = 2x U_k(x) - U_{k-1}(x)~,\mbox{ for all }k\geq 0~,
\end{equation}
and are orthogonal with respect to the Wigner semicircle measure $d
\sigma(x) = \frac{2}{\pi}\sqrt{1-x^2}\mathbf{1}_{[-1,1]}(x)d x$.

Let $A=A(G)$ denote the adjacency matrix of $G$, and define the
$n\times n$ matrix $A^{(k)}$ for $k \geq 1$:
$$ A_{u,v}^{(k)} = |\widetilde{\mathcal{W}}_{u,v}^{(k)}|\mbox{ for all }u,v\in V~.$$
That is, the entry of $A^{(k)}$ at indices $u,v$ is equal to the
number of non-backtracking walks of length $k$ from $u$ to $v$. By
definition, the matrices $A^{(k)}$ satisfy the following recurrence
relation:
\begin{equation}
  \label{eq-A-k-recursion}
  \left\{
  \begin{array}
    {l}    A^{(1)}=A~,~A^{(2)}=A^2 - d I~,\\
A^{(k+1)}= A A^{(k)} - (d-1)A^{(k-1)} \mbox{ for
    }k=2,3,\ldots~.
  \end{array}
  \right.
\end{equation}
where the last term above, $(d-1)A^{(k-1)}$, eliminates the walks
which backtrack in the $k+1$ step. We claim that:
\begin{equation}
  \label{eq-A-k-Pk-poly}
  A^{(k)}= \sqrt{d (d-1)^{k-1}}
  q_k\left(\frac{A}{2\sqrt{d-1}} \right)~\mbox{ for all }k \geq 1~,
\end{equation}
where:
\begin{equation}
  \label{eq-Pk-poly}
  q_k(x) =\sqrt{\frac{d-1}{d}} U_k(x) -
  \frac{1}{\sqrt{d(d-1)}}U_{k-2}(x)~\mbox{ for all }k\geq 1~.
\end{equation}
To see this, let $f(A,k)=\sqrt{d (d-1)^{k-1} }
  q_k\left(A/(2\sqrt{d-1})\right) $ denote the right hand side of
\eqref{eq-A-k-Pk-poly}. Substituting the polynomials $U_{-1}(x)=0$,
$U_0(x)=1$, $U_1(x)=2x$ and $U_2(x)=4x^2-1$ in \eqref{eq-Pk-poly}
implies that $f(A,1)=A=A^{(1)}$ and that $f(A,2)=A^2 - d I =
A^{(2)}$, confirming \eqref{eq-A-k-Pk-poly} for $k=1,2$. In order to
verify that \eqref{eq-A-k-Pk-poly} holds for all $k \geq 3$, recall
that $q_k(x)$ is a linear combination of the polynomials $U_{k-2}$
and $U_k$, hence it satisfies the recurrence
\eqref{eq-cheb-recursion}:
$$ q_{k+1}(x)=2x q_k(x)-q_{k-1}(x)~\mbox{ for all }k\geq 2.$$
Therefore, by induction, the following holds for all $k \geq 2$:
\begin{align}
f(A,k+1)&=\sqrt{d (d-1)^{k} }
q_{k+1}\left(\frac{A}{2\sqrt{d-1}}\right) \nonumber \\
&= \sqrt{d (d-1)^{k} } \left[\frac{A}{\sqrt{d-1}}
q_k \left(\frac{A}{2\sqrt{d-1}}\right)
- q_{k-1}\left(\frac{A}{2\sqrt{d-1}}\right)\right] \nonumber \\
&= A A^{(k)} - (d-1)A^{(k-1)} = A^{(k+1)}~,\nonumber
\end{align}
where the last inequality is by \eqref{eq-A-k-recursion}.
\begin{remark}
One can verify that the polynomials $q_k(x)$ are orthogonal
polynomials with respect to the Kesten-McKay measure
$\displaystyle{d\sigma(x) =
\frac{2d(d-1)}{\pi}\frac{\sqrt{1-x^2}}{d^2-4(d-1)x^2}}\mathbf{1}_{[-1,1]}(x)dx$.\end{remark}
Take $k \geq 1$, and recall that $A^{(k)}_{u,v}$ is the number of
non-backtracking walks of length $k$ from $u$ to $v$. Normalizing
the matrix $A^{(k)}$ as follows:
\begin{equation}\label{eq-Pk-Ak-relation}\widetilde{P}^{(k)} = \frac{A^{(k)}}{d(d-1)^{k-1}}~,
\end{equation} we obtain that
$\widetilde{P}^{(k)}$ is precisely the transition probability matrix
of a non-backtracking random walk of length $k$. Let $\mu_1=1,
\mu_2,\ldots,\mu_n$ denote the eigenvalues of $\widetilde{P}^{(k)}$,
and let \begin{equation}\label{eq-mu-def}\mu = \mu(k) =
\max\{|\mu_2|,\ldots,|\mu_n|\}~.\end{equation}
\begin{claim}\label{clm-eigenvalue-mixing}
  Let $\widetilde{P}_{ij}^{(k)}$ and $\mu(k)$ be as above. The following holds:
\begin{equation}\label{eq-Pij-mu-relation}
\frac{\mu(k)}{n} \leq \max_{i,j}\Big|\widetilde{P}^{(k)}_{i j} - \frac{1}{n} \Big| \leq
\mu(k)~.\end{equation}
\end{claim}
\begin{proof}
The vector $v_1 = \frac{1}{\sqrt{n}}(1,\ldots,1)$ is an
eigenvector of $\widetilde{P}^{(k)}$ corresponding to its largest eigenvalue $\mu_1=1$, and
therefore:
\[ \max_{i,j} \left| \widetilde{P}^{(k)}_{ij} - \frac{1}{n} \right|
    = \max_{i,j} \left| \left\langle \left( \widetilde{P}^{(k)} - v_1 \otimes v_1
    \right) e_i, e_j \right\rangle \right|
    \leq \max_{|u|=|v|=1} \left| \left\langle \left( \widetilde{P}^{(k)}- v_1 \otimes v_1\right) u, v \right\rangle \right| = \mu(k)
    ~.  \]
On the other hand:
\[ \max_{i,j} \left| \widetilde{P}^{(k)}_{ij} - \frac{1}{n} \right|
    \geq \frac{1}{n}
        \sqrt{\sum_{i,j} \big| \widetilde{P}^{(k)}_{ij} - \frac{1}{n} \big|^2}
    = \frac{1}{n} \sqrt{\sum_{2\leq s \leq n} \mu_s^2} \geq \mu(k)/n \, \text{.} \]
\end{proof}
We deduce that:
\begin{equation}
  \label{eq-rho-equals-second-ev}
  \widetilde{\rho} = \limsup_{k\to\infty} \mu(k)^{1/k} = \max_{2\leq i \leq n} \limsup_{k\to\infty} |\mu_i(k)|^{1/k}~,
\end{equation}
and it remains to compute the right hand side above. By \eqref{eq-A-k-Pk-poly} and \eqref{eq-Pk-Ak-relation}, the following holds
for all $i \in [n]$:
$$ \mu_i = \frac{1}{\sqrt{d(d-1)^{k-1}}} q_k \Big(\frac{\lambda_i}{2\sqrt{d-1}}\Big) ~,$$
where $\lambda_i$ are the eigenvalues of $A$.
Therefore, the proof of the theorem will follow from the next lemma:
\begin{lemma}
The polynomials $q_k$, defined in \eqref{eq-Pk-poly}, satisfy:
  \[ \limsup_{k \to \infty} |q_k(x)|^{1/k} = \psi(|x|) =
\begin{cases}
1,                  & -1 \leq x \leq 1~, \\
|x| + \sqrt{x^2-1},   & x \in \mathbb{R}\setminus[-1,\,1]~.
\end{cases} \]
\end{lemma}
\begin{proof}
If $x \in [-1, 1]$, then $x = \cos \theta$ for some $\theta \in [0,\pi]$, and hence:
\begin{equation}\label{eq-qk-x-in-[-1,1]} q_k(x) = \sqrt\frac{d-1}{d} \frac{\sin ((k+1)\,\theta)}{\sin \theta}
    - \frac{1}{\sqrt{d(d-1)}} \frac{\sin ((k-1)\,\theta)}{\sin \theta}
    \, ~.\end{equation}
Therefore:
\[ |q_k(x)| \leq \sqrt\frac{d-1}{d} (k+1) + \frac{1}{\sqrt{d(d-1)}} (k-1) \]
and $\limsup_{k\to\infty} |q_k(x)|^{1/k} \leq 1$. The reverse inequality follows
from an appropriate subsequence $k_j$ for which the right hand side of \eqref{eq-qk-x-in-[-1,1]}
is bounded from below by some $c=c(\theta) > 0$.

It remains to treat $x \notin [-1,1]$. In this case, $x = (z+z^{-1})/2$ for $z = x + \sign(x) \sqrt{x^2-1} \notin [-1,1]$.
Setting $z = \sign(x)\mathrm{e}^{\theta}$ for some real $\theta$, we get $x = \sign(x)\cos(i\theta)$, and therefore:
\begin{align} q_k(x) &= \sign(x)^k\left(\sqrt\frac{d-1}{d}
\frac{\sin ((k+1)i\theta)}{\sin (i \theta)}
    - \frac{1}{\sqrt{d(d-1)}} \frac{\sin ((k-1)i\theta)}{\sin
    (i\theta)}\right)
    \nonumber \\
    &=
\sqrt\frac{d-1}{d} \frac{z^{k+1}-z^{-(k+1)}}{z-z^{-1}} -
    \frac{1}{\sqrt{d(d-1)}} \frac{z^{k-1}-z^{-(k-1)}}{z-z^{-1}} ~,\nonumber \end{align}
and $\limsup |q_k(x)|^{1/k} = \lim |q_k(x)|^{1/k} = |z|$.

This completes the proof of the lemma and of Theorem \ref{thm-1}.
\end{proof}

\begin{proof}[{\em \textbf{Proof of Corollary \ref{cor-mixing-ratio}}}]
Let $\lambda$ denote the largest absolute value of a nontrivial
eigenvalue of $G$.
Note that $\psi(x)$, as defined in Theorem \ref{thm-1}, satisfies the following properties:
\begin{equation} \label{eq-psi-props}
\left\{\begin{array}{l} \psi \mbox{ is strictly monotone increasing
on
}[1,\infty)~,~\psi(1)=1~,~\displaystyle{\frac{\psi(x)}{x}\mathop{\longrightarrow}_{x\to\infty}
2}~,\\
 \displaystyle{\frac{\psi\left(\frac{x}{2\sqrt{d-1}}\right)}{\sqrt{d-1}} =
 \frac{x + \sqrt{x^2-4d+4}}{2(d-1)} \leq \frac{x}{d}}~\mbox{ for all
 $d$ and all $2\sqrt{d-1}\leq x \leq d$}
 ~,\\
 \noalign{\medskip}
 \psi\left(\frac{d}{2\sqrt{d-1}}\right)= \sqrt{d-1}~.
\end{array}\right.
\end{equation}
Therefore, if $\lambda > 2\sqrt{d-1}$, Theorem \ref{thm-1}
implies that
$\widetilde{\rho} = \psi\left(\frac{\lambda}{
2\sqrt{d-1}}\right)/\sqrt{d-1}$, and that:
$$\frac{\lambda}{2(d-1)} <
\widetilde{\rho} \leq \frac{\lambda}{d}~.$$ As $\rho=\lambda / d$,
we obtain \eqref{eq-mixing-rates-ratio}. Furthermore, as $\lambda$
decreases to $2\sqrt{d-1}$,
$\psi\left(\frac{\lambda}{2\sqrt{d-1}}\right)$ tends to $1$, implying
that $\widetilde{\rho} \to \frac{1}{\sqrt{d-1}}$, and
$\widetilde{\rho}/\rho \to \frac{d}{2(d-1)}$.

It remains to handle the case $\lambda \leq 2\sqrt{d-1}$. To this
end, recall the following result of Nilli \cite{Nilli}, which
implies the Alon-Boppana Theorem:
\begin{theorem}
[\cite{Nilli}] If $G$ is a simple undirected $d$-regular graph with
diameter at least $2(k+1)$, then the second largest eigenvalue of
$G$, $\lambda_2$, satisfies $\lambda_2 \geq
2\sqrt{d-1}-\frac{2\sqrt{d-1}-1}{k+1}$.
\end{theorem}
As the diameter of a $d$-regular graph on $n$ vertices is at least
$(1-o(1))\log_{d-1} n$, we deduce that in the above case, if
$d=n^{o(1)}$ then $\lambda = (1-o(1))2\sqrt{d-1}$. In this case, by
Theorem \ref{thm-1} we have $\widetilde{\rho}=1/\sqrt{d-1}$, and
$\widetilde{\rho}/\rho=\frac{d}{2(d-1)}+o(1)$.
\end{proof}

\begin{remark}
  Examining the trace of the square of the adjacency matrix of a graph, it is easy
  to see that for every $d$-regular graph on $n$ vertices, the second largest eigenvalue in absolute
  value is at least $\sqrt{\frac{d(n-d)}{(n-1)}}$. It thus follows that if $d=o(n)$ then $\widetilde{\rho}\leq(1+o(1))\rho$.
\end{remark}
\begin{remark}
  For $d$-regular graphs with $d=\Theta(n)$ the mixing rate of the simple random walk may indeed
  be faster than that of the non-backtracking random walk. For instance, if $G$ is the complete
  graph on $n$ vertices, $K_n$, then by Theorem \ref{thm-1}, $\widetilde{\rho}=\frac{1}{\sqrt{n-2}}$, and
  $\rho = \frac{1}{n-1}$.
\end{remark}

\section{Random walks and the balls and bins paradigm}\label{sec::balls-bins}
\noindent \textbf{Proof of Theorem \ref{thm-2}:} Let $G$ be as
described in Theorem \ref{thm-2}. The following definition of the
 mixing-time of a non-backtracking random walk on $G$ corresponds
  to an $L_\infty$ distance of $1/n^2$ between
$\pi$ and $\widetilde{P}^{(k)}_u$, for all $u\in V$:
\begin{equation}\label{eq-fine-mixing-time-def}
\tau  = \min_t \left\{ \Big|\widetilde{P}^{(k)}_{u v} - \frac{1}{n}
\Big| \leq \frac{1}{n^2} ~\mbox{ for all $u,v\in V$ and $k \geq
t$}\right\}~.
\end{equation}
Theorem \ref{thm-1} implies that a non-backtracking random walk on
$G$ converges to the uniform distribution at a mixing-rate of
$\widetilde{\rho}=\psi\left(\frac{\lambda}{2\sqrt{d-1}}\right)/\sqrt{d-1}$,
and we deduce that $\tau = O(\log n)$ (by usual arguments linking
the mixing-rate to the mixing-time).

The proof of Theorem \ref{thm-2} will follow from the next two
lemmas, which we prove using first and second moment arguments (see,
e.g., \cite{ProbMethod}), combined with some additional ideas.
\begin{lemma}
\label{lem-balls-upper-bound} Let $G$ be as in Theorem \ref{thm-2}.
With high probability, a non-backtracking random walk of length $n$
on $G$ does not visit a vertex more than $(1+o(1))\frac{\log
n}{\log\log n}$ times.
\end{lemma}

\begin{lemma}\label{lem-balls-lower-bound} Let $G$ be as in Theorem \ref{thm-2}.
With high probability, a non-backtracking random walk of length $n$
on $G$ visits some vertex at least $\left(1+o(1)\right)\frac{\log
n}{\log\log n}$ times.
\end{lemma}

The key element in the proofs of both lemmas is showing that the
number of times that a non-backtracking random walk visits some
vertex, or some pair of vertices, is governed by visits at locations
which are at least $\tau$ apart. This implies a behavior which is
essentially the same as the one in the balls and bins experiment.

\begin{proof}[{\em \textbf{Proof of Lemma
\ref{lem-balls-upper-bound}}}] Let $u,v \in V$ denote two vertices,
so that either $u=v$ or the distance between $u$ and $v$ in $G$ is
at least \begin{equation}\label{eq-L-def}L = 10\log_{d-1}\log
n~,\end{equation} and let $\widetilde{P}^{(\ell)}_{u v}$ denote the
probability that a non-backtracking random walk of length $\ell$ on
$G$, which starts at $u$, ends in $v$. We claim that:
\begin{equation} \label{eq-u-v-path-bound} \widetilde{P}^{(\ell)}_{u v} \leq
\left\{\begin{array} {ll} (d-1)/(\log n)^5 & \mbox{If }\ell
< \tau~,\\
(1+n^{-1})/n&\mbox{If }\ell \geq \tau~.
\end{array}\right.
\end{equation}
The case $\ell \geq \tau$ follows directly from the definition
\eqref{eq-fine-mixing-time-def} of the mixing time $\tau$. For the
case $\ell < \tau$, let $W=(u=w_0,w_1,\ldots,w_\ell)$ denote a
non-backtracking random walk of length $\ell$ on $G$, starting at
$u$. The choice of $u,v$ and the fact that $g$, the girth of $G$, is
at least $L$ (this applies to the case $u=v$), imply that there is
no non-empty path between $u,v$ of length shorter than $L$.
Therefore, if $\ell < L$ then $\Pr[w_\ell = v] = 0$. Otherwise, let
$h=\lfloor \frac{L-1}{2}\rfloor$, and notice that the neighborhood
of $v$ up to distance $h$ is precisely a $d$-regular tree (as $L
\leq g$). Let $U$ denote the $d(d-1)^{h-1}$ leaves of this tree.
Since the random walk $W$ cannot backtrack, the event $w_\ell=v$
implies that $w_{\ell-h}\in U$, hence:
\begin{align}\widetilde{P}_{u v}^{(\ell)} = \Pr[w_\ell = v] &= \Pr[w_\ell = v ~|~
w_{\ell-h} \in U] \Pr[w_{\ell-h}\in U] \nonumber \\
&\leq \Pr[w_\ell = v ~|~ w_{\ell-h} \in U] = (d-1)^{-h} \leq
\frac{d-1}{(\log n)^5}~.\nonumber
\end{align}
Let $\epsilon > 0$, and set $k = (1+\epsilon)\frac{\log n}{\log\log
n}$. Consider a non-backtracking random walk of length $n$ on $G$,
$W=(w_0,w_1,\ldots,w_n)$, where $w_0$ is a fixed vertex of $V$. For
each vertex $v \in V$, and for each $t\in\{0,\ldots,k\}$, define the
following event:
$$A_{v,t} = \left(\begin{array}{l}
\mbox{$W$ visits $v$ at least $k$ times at some indices }
1\leq i_1<\ldots<i_k<\ldots~,\\
\mbox{and }|\{j\in[k-1]: i_{j+1}-i_j < \tau\}| = t~.
\end{array}
\right) ~.$$ That is, $A_{v,t}$ describes the event in which
precisely $t$ of the first $k-1$ segments of $W$, which are bounded
by consecutive visits to $v$, are of length smaller than $\tau$.
Considering all the possible ways to choose indices $i_0,\ldots,i_k$
according to the definition of $A_{v,t}$, we derive the following
from \eqref{eq-u-v-path-bound}: \begin{equation}
\label{eq-A-vt-upper-bound} \Pr[A_{v,t}] \leq
\binom{n}{k-t}\binom{k-1}{t}\tau^{t}
\left(\frac{1+n^{-1}}{n}\right)^{k-t} \left(\frac{d-1}{(\log
n)^5}\right)^{t}~.\end{equation} For $0 \leq t < k-1$, replacing $t$
by $t+1$ in the right hand side of \eqref{eq-A-vt-upper-bound}
results in a multiplicative factor of:
$$ \frac{(k-t)(k-t-1)}{(n-k+t+1)(t+1)}\cdot \frac{\tau n}{1+n^{-1}}\cdot\frac{d-1}{(\log n)^5}
= O\left( \frac{k^2\tau}{(\log n)^5}\right) = o(1)~.$$ Therefore,
the largest term is obtained for $t=0$. Letting $A_v =
\cup_{t=0}^{k-1}A_{v,t}$ denote the event that $W$ visits the
vertex $v$ at least $k$ times, we get:
$$\Pr[A_v ] \leq k \binom{n}{k}
\left(\frac{2}{n}\right)^{k} \leq \frac{2^k}{(k-1)!}=o(1/n)~,$$
where the last inequality is by the assumption on $k$. Therefore,
$\Pr[\cup_{v\in V}A_v]=o(1)$, and with high probability, $W$ does
not visit any vertex of $V$ more than $k$ times.
\end{proof}
\begin{proof}[{\em \textbf{Proof of Lemma
\ref{lem-balls-lower-bound}}}] Let $\epsilon > 0$, and set $k =
\lceil (1-\epsilon)\frac{\log n}{\log\log n}\rceil$. Let $W$ denote
a non-backtracking random walk of length $n$ on $G$,
$W=(w_0,w_1,\ldots,w_n)$, where $w_0$ is a fixed vertex of $V$. We
wish to show that, with high probability, $W$ visits some vertex $v
\in V$ at least $k$ times. We will show that, in fact, this
statement holds even if we restrict ourselves to a predefined subset
of the vertices $U \subset V$, and in addition, restrict the pattern
of the visiting locations.

Let $U \subset V$ denote a set of vertices of $G$ of size
\begin{equation}
  \label{eq-U-size-lower-bound}
  |U| = \lceil n / \left(d(\log n)^{10}\right) \rceil~,
\end{equation}
so that the distance between any pair of vertices $u,v\in U$ is at
least $L=10\log_{d-1}\log n$ (as defined in \eqref{eq-L-def}). To
see that such a set $U$ indeed exists, notice that the number of
vertices, whose distance from some $u \in U$ is at most $L-1$, does
not exceed $\sum_{i=0}^{L-1}d(d-1)^i \leq d(d-1)^L$. Therefore, a
greedy algorithm which begins with an empty set, and repeatedly adds
a new legal vertex to $U$, always succeeds in producing a set of
size at least $n / \left(d(\log n)^{10}\right)$.

The restriction we impose on the pattern of visits is defined next:
\begin{definition}Let $T\subset \binom{[n]}{k}$ denote a set of $k$ indices in
$[n]$. We say that $T$ is a {\bf $k$-pattern} iff $T \cap [2\tau] =
\emptyset$ and $|i-j| > 2\tau$ for all $i,j \in T$. In other words,
the value of the elements of $T$, and the pairwise distances between
these elements, all exceed $2\tau$.\end{definition} The above
definition implies that, if $T$ is a $k$-pattern, then for all
$i\in[n]$, there is at most one element $j\in T$ so that $|i-j|\leq
\tau$. This makes it useful to define the correlation between
$k$-patterns as follows:
\begin{definition}
  Let $T_1$ and $T_2$ denote two $k$-patterns. The correlation between $T_1$ and $T_2$,
  $\delta(T_1,T_2)$, is defined as the number of pairs in
  $T_1\times T_2$ with distance at most $\tau$:
  $$\delta(T_1,T_2)=\left|\big\{(a,b) \in T_1 \times T_2: |a-b|\leq
  \tau\big\}\right|~.$$
\end{definition}
Let $\mathcal{K}$ denote the collection of all $k$-patterns, and
notice that:
\begin{equation}\label{eq-K-size}
|\mathcal{K}| = \binom{n-2\tau k}{k}~.
\end{equation}
Define the following set of indicator variables for all $u \in U$
and $T\in\mathcal{K}$:
\begin{equation}
  \label{eq-a-u,T-def}
X_{u,T}=\left\{\begin{array}{ll}1 & \mbox{If $w_i = u$ for all $i
\in T$}~,\\
0&\mbox{otherwise.}\end{array}\right.
\end{equation}
In other words, $X_{u,T}$ is the indicator for the event according
to which the non-backtracking walk $W$ visits $u$ in all the
time-points specified by $T$. By definition, the first of these
time-points exceeds $\tau$, and the same holds for the distance
between each consecutive pair of these time-points, and by the
definition of $\tau$ we deduce that:
\begin{equation}
  \label{eq-pr-X-uT}
\left(\frac{1-n^{-1}}{n}\right)^k \leq \Pr[X_{u,T}=1] \leq
\left(\frac{1+n^{-1}}{n}\right)^k~.
\end{equation}
Setting $X = \sum_{u,T}X_{u,T}$, we get:
\begin{equation}\label{eq-X-exp-lower-bound}
\mathbb{E}X \geq |U| \binom{n-2\tau k}{k}
\left(\frac{1-n^{-1}}{n}\right)^k = n^{\epsilon-o(1)}~,
\end{equation}
where the last equality is by the definition of $k$ and
\eqref{eq-U-size-lower-bound}.

In order to show that $X$ is concentrated around its expected value,
we consider its second moment. Let $u,v \in U$ so that $u \neq v$,
and let $t \in \{0,1,\ldots,k\}$. Take $T_1,T_2 \in \mathcal{K}$ so
that $\delta(T_1,T_2)=t$. By the definition of $U$, the distance
between $u$ and $v$ is at least $L$. Hence, if $|a-b|< L$ for some
$(a,b)\in T_1 \times T_2$, then the events $(X_{u,T_1}=1)$ and
$(X_{v,T_2}=1)$ are disjoint. Otherwise, consider the probability of
the event $(X_{u,T_1}=1) \wedge (X_{v,T_2}=1)$. By
\eqref{eq-u-v-path-bound}, the largest of each of the $t$ pairs of
indices $(a_i,b_i) \in T_1 \times T_2$, which satisfy $|a_i-b_i|\leq
\tau$, contributes a probability of at most $(d-1)/(\log n)^5$ to
this event. The definition of $\tau$ implies that each of the
remaining indices contributes a probability of at most
$(1+n^{-1})/n$ for visiting the required vertex (either $u$ or $v$),
and altogether:
\begin{equation}\label{eq-pr-wedge-X-uT}\Pr[X_{u,T_1}=1
\wedge X_{v,T_2}=1] \leq \left(\frac{1+n^{-1}}{n}\right)^{2k-t}
 \left(\frac{d-1}{(\log n)^5}\right)^t~.\end{equation} Combining \eqref{eq-pr-X-uT} and
\eqref{eq-pr-wedge-X-uT} gives:
\begin{align}
&\sum_{T_1\in \mathcal{K}} \mathop{\sum_{T_2\in
\mathcal{K}}}_{\delta(T_1,T_2)=t} \cov(X_{u,T_1},X_{v,T_2})
\nonumber \\ &\leq \ \binom{n-2\tau k}{k} \binom{k}{t} (2\tau)^t
\binom{n-2\tau k}{k-t} \left(\left(\frac{1+n^{-1}}{n}\right)^{2k-t}
 \left(\frac{d-1}{(\log n)^5}\right)^t - \left(\frac{1-n^{-1}}{n}\right)^{2k}
\right)\nonumber \\
&= \binom{n-2\tau k}{k} \binom{k}{t} (2\tau)^t \binom{n-2\tau
k}{k-t} n^{-2k} \left(\left(1+n^{-1}\right)^{2k}
 \left(\frac{(1+o(1))(d-1)n}{(\log n)^5}\right)^t - \left(1-n^{-1}\right)^{2k}
\right)~.\nonumber
\end{align}
Let $C_{u v}(t)$ denote the right hand side in the above inequality.
Since $(1+n^{-1})^{2k}$ and $(1-n^{-1})^{2k}$ both tend to $1$ as
$n$, and hence $k$, tend to $\infty$, the following holds for all $t
\geq 1$: \begin{align}\frac{C_{u v}(t+1)}{C_{u v}(t)} &=
\frac{(k-t)^2(2\tau)}{(t+1)(n-(2\tau+1) k + t+1)}\cdot (1+o(1))
\frac{(d-1)n}{(\log n)^5} = O\left(\frac{k^2 \tau}{(\log
n)^5}\right) = o(1)~. \nonumber
\end{align}
In particular, for a sufficiently large $n$ we deduce that
\begin{equation}\label{eq-Cuv-sum}\sum_{t =1}^k C_{u v}(t) \leq 2C_{u v}(1)~,\end{equation} and it remains to
examine $C_{u v}(t)$ for $t\in\{0,1\}$:
\begin{align}\label{eq-Cuv(0)-bound}C_{u v}(0) &= \binom{n-2\tau k}{k}^2
\left(\left(\frac{1+n^{-1}}{n}\right)^{2k} -
\left(\frac{1-n^{-1}}{n}\right)^{2k} \right)\nonumber \\
&\leq \left(\frac{\mathbb{E}X}{|U|}\right)^2
\left(\frac{n}{1-n^{-1}}\right)^{2k} \cdot\frac{2}{n^2}\cdot2k
\left(\frac{1+n^{-1}}{n}\right)^{2k-1} = O\left(\frac{k}{n}
\left(\frac{\mathbb{E}X}{|U|}\right)^2\right)~,\\
\label{eq-Cuv(1)-bound}C_{u v}(1) &= \binom{n-2\tau
k}{k} \cdot 2 \tau k  \binom{n-2\tau k}{k-1}
n^{-2k}\cdot(1+o(1))\frac{(d-1)n}{(\log n)^5}
\nonumber \\
&\leq \left(\frac{\mathbb{E}X}{|U|}\right)^2
\left(\frac{n}{1-n^{-1}}\right)^{2k} \cdot \frac{2\tau
k}{n^{2k}}\cdot \frac{(1+o(1))k(d-1)}{(\log n)^5} =
O\left(\frac{\tau k^2}{(\log n)^5}
\left(\frac{\mathbb{E}X}{|U|}\right)^2\right)~.
\end{align}
By \eqref{eq-Cuv-sum}, \eqref{eq-Cuv(0)-bound} and
\eqref{eq-Cuv(1)-bound} we get:
\begin{equation}\label{eq-cov-u-neq-v} \sum_{u\in U}
\mathop{\sum_{v\in U}}_{u \neq v} \sum_{(T_1,T_2)\in\mathcal{K}^2}
\cov(X_{u,T_1},X_{v,T_2}) \leq \sum_{u\in U} \mathop{\sum_{v\in
U}}_{u \neq v} \sum_{t=0}^k C_{u v}(t) =
o(\left(\mathbb{E}X\right)^2)~.\end{equation} Next, take $u \in U$,
and consider all $k$-patterns $T_1 \neq T_2$ which contain $l\in
\{0,\ldots,k-1\}$ common indices, and whose correlation,
$\delta(T_1,T_2)$, is some $t \in \{l,\ldots,k\}$. The following
holds:
\begin{align} & \sum_{T_1 \in \mathcal{K}} \mathop{\mathop{\sum_{T_2 \in
\mathcal{K}}}_{|T_1 \cap T_2|=l}}_{\delta(T_1,T_2)=t}
\cov(X_{u,T_1},X_{u,T_2})  \nonumber \\
& \leq \binom{n-2\tau k}{k} \binom{k}{t} \binom{t}{l} (2\tau)^{t-l}
\binom{n-2\tau k}{k-t} \left(\left(\frac{1+n^{-1}}{n}\right)^{2k-t}
 \left(\frac{d-1}{(\log n)^5}\right)^{t-l} - \left(\frac{1-n^{-1}}{n}\right)^{2k}
\right)\nonumber \\
& \leq (1+o(1))\binom{n-2\tau k}{k} \binom{k}{t} \binom{t}{l}
\binom{n-2\tau k}{k-t} n^{-2k+t}
 \left(\frac{2\tau(d-1)}{(\log n)^5}\right)^{t-l}~.\label{eq-Cu-exp}
\end{align}
Let $C_u(l,t)$ denote the final expression of \eqref{eq-Cu-exp}. For
all $l$ and $t$ so that $l\leq t < k$ we have:
\begin{equation} \frac{C_u(l,t+1)}{C_u(l,t)} = \frac{(1+o(1))(k-t)^2}{t-l+1} \cdot
\frac{2\tau(d-1)}{(\log n)^5} = O\left(\frac{k^2 \tau}{(\log
n)^5}\right) = o(1)~.\label{eq-Cu-t-increment} \end{equation} This
implies that the leading order term in the sum $\sum_{t=l}^k
C_u(l,t)$ is $C_u(l,l)$. Next,
$$ \frac{C_u(l+1,l+1)}{C_u(l,l)} = \frac{(1+o(1))(k-l)^2}{l+1}~,$$
hence, if we define: $$l_0 = k - 2\sqrt{k}~,~l_1 = k
-\frac{1}{2}\sqrt{k}~,$$ then the following holds:
\begin{equation} \left\{ \begin{array}
  {ll}
  \frac{C_u(l+1,l+1)}{C_u(l,l)} \geq 4+o(1) & \mbox{if }l \leq l_0~,\\
    \noalign{\medskip}
\frac{C_u(l+1,l+1)}{C_u(l,l)} \leq \frac{1}{4}+o(1) & \mbox{if }l
\geq l_1~.
\end{array}\right.\label{eq-Cu-l-increment} \end{equation} On the
other hand, for every $l \in [l_0,l_1]$ we have:
\begin{align}C_u(l,l)&= (1+o(1))\binom{n-2\tau k}{k}\binom{k}{l} \binom{n-2\tau
k}{k-l}n^{-2k+l} \nonumber \\
&\leq (1+o(1))\frac{\mathbb{E}X}{|U|}\left(\frac{\mathrm{e}^2 k
(n-2\tau k)}{(k-l)^2 n}\right)^{k-l} \leq
(1+o(1))\frac{\mathbb{E}X}{|U|}\left(4\mathrm{e}^2+o(1)
\right)^{2\sqrt{k}}\nonumber \\
& = (1+o(1))\frac{\mathbb{E}X}{|U|}n^{o(1)} =
o\left(\frac{(\mathbb{E}X)^2}{|U|n^{\epsilon/2}}\right)
\label{eq-Cu-l-upper-bound}~,\end{align} where the last equality is
by \eqref{eq-X-exp-lower-bound}.
 We
deduce from \eqref{eq-Cu-t-increment}, \eqref{eq-Cu-l-increment} and
\eqref{eq-Cu-l-upper-bound} that for all sufficiently large values
of $n$:
\begin{align}\sum_{l=0}^{k-1} \sum_{t=l}^k C_u(l,t) &\leq
2\sum_{l=0}^{k-1} C_u(l,l) \leq 4 C_u(l_0,l_0) + 4C_u(l_1,l_1) +
2\sum_{l=l_0}^{l_1}C_u(l,l)\nonumber \\
&= o\left(\sqrt{k}\frac{(\mathbb{E}X)^2}{|U|n^{\epsilon/2}}\right) =
o\left(\frac{(\mathbb{E}X)^2}{|U|}\right) ~,\nonumber\end{align} and
thus:
\begin{align} &\sum_{u \in U} \sum_{T_1 \in \mathcal{K}}
\mathop{\sum_{T_2\in\mathcal{K}}}_{T_1 \neq T_2}
\cov(X_{u,T_1},X_{u,T_2}) \leq \sum_{u\in U} \sum_{l=0}^{k-1}
\sum_{t=l}^k C_u(l,t) = o\left((\mathbb{E}X)^2\right)~.
\label{eq-cov-u}
\end{align}
Combining \eqref{eq-cov-u-neq-v} and \eqref{eq-cov-u} (and recalling
that $\mathbb{E}X = \omega(1)$) gives:
$$\var(X) \leq \mathbb{E}X + \sum_{u,T_1}
\mathop{\sum_{v,T_2}}_{(u,T_1)\neq (v,T_2)}
\cov(X_{u,T_1},X_{v,T_2}) = o((\mathbb{E}X)^2)~,$$ and Chebyshev's
inequality implies that:
$$ \Pr[X = 0] \leq \frac{\var(X)}{(\mathbb{E}X)^2} = o(1)~.$$ This
completes the proof of Lemma \ref{lem-balls-lower-bound} and of
Theorem \ref{thm-2}.
\end{proof}

We note that the $\Omega(\log\log n)$ requirement on the girth of
$G$ in Theorem \ref{thm-2} is tight, as there are
$(n,d,\lambda)$-graphs with girth $g$, where a non-backtracking
random walk visits some vertex at least $\Omega(\frac{\log n}{g})$
times almost surely. This is stated in the next claim.
\begin{claim}\label{clm-low-girth}
Let $G$ be a $d$-regular graph on $n$ vertices, in which each vertex
is contained in a cycle of length $g=g(n)$. If $k=k(n)$ satisfies:
$$ k = \frac{\log_{d-1} \left(n/\log n\right) - \omega(1)}{g}~,$$ then,
with high probability, a non-backtracking random walk of length $n$
on $G$ visits some vertex at least $k$ times. In particular, such a
walk almost surely visits some vertex $\Omega(\frac{\log n}{g})$
times.
\end{claim}
\begin{proof}
For each $v \in V$, let $C_v$ denote a cycle of length $g$ which
contains $v$ in $G$. Let $W=(w_0,w_1,\ldots,w_n)$ denote a
non-backtracking random walk of length $n$ on $G$, and divide $W$
into $T=\lfloor n/kg \rfloor$ disjoint segments, $I_1,\ldots,I_T$,
each of length $kg$:
$$I_j = \left(w_{(j-1)k g},\ldots,w_{j k
g - 1}\right)~\mbox{for all }j \in T~.$$ Define the following event
for each $j \in [T]$:
$$ A_j = \left(\begin{array}{l} \mbox{The segment $I_j$ of $W$ is precisely $k$ consecutive walks}\\
\mbox{along the same cycle $C_v$, where $v = w_{(j-1)kg}$.}
\end{array}\right)~.
$$
To prove the claim, it suffices to show that, with high probability,
at least one of the events $A_j$ ($j\in[T])$ occurs. Since these
events are independent, and $\Pr[A_j] = (d-1)^{-k g}$ for all $j$,
we get:
$$\Pr[\cap_{j=1}^{T}A_j^{c}] = \left(1- (d-1)^{-k g}\right)^T \leq
\exp(-T(d-1)^{-k g})~.$$ The choice of $k$ ensures that $ T(d-1)^{-k
g} = \omega(1)$, and the result follows.
\end{proof}

\begin{remark}
Theorem \ref{thm-2} stated that the maximal load in a
non-backtracking walk of length $n$ on a $d$-regular expander of
high girth is $(1+o(1))\frac{\log n}{\log\log n}$ with high probability,
similar to the
maximal load in the classical balls and bins experiment. In contrast to this,
a simple calculation shows that a typical simple random walk of length $n$,
on any $d$-regular graph for a fixed $d$, has a maximal load of $\Omega(\log n)$.
This can be seen as a special case of Claim \ref{clm-low-girth},
taking $g=2$: the probability that the simple random walk traverses
the same edge repeatedly for, say, $M=\frac{1}{2} \log_d n$
consecutive steps, is $1/\sqrt{n}$. Dividing the walk to disjoint
segments of length $M$ implies that, with probability $1-o(1)$, at
least one segment exhibits this behavior, thus the maximal load is
at least $M$.
\end{remark}

\begin{remark} The classical Birthday Paradox states that, when
throwing balls to $n$ bins, independently and uniformly at
random, we expect a collision after $\Theta(\sqrt{n})$ balls (see,
e.g., \cite{Feller}). Relating this to random walks on expanders,
one may ask when do simple and non-backtracking random walks on
expanders self-intersect. Clearly, most {\em simple} random
walks on an expander encounter a collision after $O(1)$ steps
(the first time at which an edge is traversed twice in a row). An
argument similar to the one used in the proof of Claim
\ref{clm-low-girth} shows that, for every small $\epsilon > 0$,
there are $(n,d,\lambda)$-graphs with girth $g = \epsilon
\log_{d-1} n$, on which a non-backtracking random walk will self
intersect after at most $n^{\epsilon+o(1)}$ steps almost surely.
Similarly, for $g=o(\log n)$, there are such graphs where the
self-intersection time of the non-backtracking random walk is at
most $(d-1)^{g+o(1)}$.
\end{remark}

\section{Concluding remarks and open
problems}\label{sec::conclusion}
\begin{itemize}
  \item We have shown that a non-backtracking random walk on every
  connected and non-bipartite $d$-regular graph $G$, where $d \geq 3$, converges to the
  uniform distribution, and computed its precise mixing-rate. We obtained that this mixing-rate
   is always asymptotically at least as fast as that of the simple random walk on the
   same graph provided $d=o(n)$ (and is faster provided $d\leq n^{o(1)}$), and their
   ratio may reach up to $2(d-1)/d$.
  \item As an application, we showed
   that if $G$ is a high-girth $d$-regular expander on $n$ vertices,
   for some fixed $d \geq 3$,
   then the maximal load while sampling $n$ consecutive positions of a
   non-backtracking random walk on $G$ is almost surely $(1+o(1))\frac{\log n}{ \log\log n}$,
   similar to the maximal load in the classical balls and bins experiment.
   Performing a simple random walk, instead of a non-backtracking one, results in a maximal
   load of $\Omega(\log n)$ with high probability.
  \item Following the Poisson approximations in the balls and bins model, it would be interesting to establish the precise
  distribution of a sample of $n$ consecutive positions of a non-backtracking random walk on an expander of
  high girth.
  \item The well known {\em power-of-two} result (\cite{ABKU}, see also \cite{MitzenmacherUpfal}, Chapter 14) states that
  if $n$ balls are thrown into $n$ bins, where each ball is placed in the
  least loaded bin, out of two independently chosen random ones, then the maximal load
  decreases from $\Theta(\frac{\log n}{\log\log n})$ to
  $\Theta(\log\log n)$. Let $W_1$ and $W_2$ denote two non-backtracking random
  walks on an expander of high girth, and suppose that in each step we
  are given a choice between the two current locations of $W_1$ and $W_2$, and
  pick the least loaded one. Does the maximal load decrease
  from $\Theta(\frac{\log n}{\log\log n})$ to $\Theta(\log\log n)$
  in this setting as-well?

  \item One way of proving the above power-of-two result in the balls and
  bins model is to consider the Erd\H{o}s-R\'enyi random graph process
  $\mathcal{G}^t$, $t\in\{0,1,\ldots,\binom{n}{2}\}$ (where $\mathcal{G}^0$ is the empty graph on $n$ vertices,
  and in each step a new edge is
  added, uniformly chosen over all missing edges; see, e.g., \cite{RandomGraphs}, Chapter
  2). Each pair of bins corresponds to a uniformly chosen edge in the graph (we
  may ignore self-loops or repeating edges, as we are dealing with a
  linear number of balls). Selecting a bin corresponds to
  choosing an orientation for this edge. One can show that the greedy
  online algorithm, which orients an edge towards the vertex with the lower
  in-degree, gives an overall maximal in-degree of $O(\log\log
  n)$ with high probability. This is based on the following properties of $\mathcal{G}^t$, which hold
  with high probability for all $t \leq \alpha n $, where $0<\alpha < \frac{1}{2}$ is a constant:
  \begin{enumerate}[(1)]
    \item \label{item-log-conn} Each connected component of $\mathcal{G}^t$ is of logarithmic size.
    \item \label{item-degen} For some fixed $\delta$, the average degree of
    every induced subgraph of $\mathcal{G}^t$ is at most $\delta$.
  \end{enumerate}
  The above discussion suggests the following approach: let $G$ be a
  $d$-regular expander of high girth, for some fixed $d \geq 3$, and let $W_1$ and $W_2$ denote
  two non-backtracking random walks on $G$. Define a random (multi) graph
  process by adding the edge $(W_1(t),W_2(t))$ at step $t$,
  where $W_i(t)$ is the position of $W_i$ at time $t$. This can be
  viewed as a certain de-randomization of the random graph process, where
  the graph at time $\Theta(n)$ is produced using only $\Theta(n)$
  random bits (instead of $\Theta(n\log n)$ bits). This model, on its own account, seems interesting, with respect to the commonly studied
  questions on graph processes, e.g., whether there exists a sharp threshold for the
  appearance of a giant component. In particular, proving that properties
  \eqref{item-log-conn} and \eqref{item-degen} hold for this graph
  process for all $t \leq \alpha n$ and some $0<\alpha < \frac{1}{2}$ will
  imply a positive answer to the previous question, regarding the power-of-two with non-backtracking random
  walks.

\end{itemize}

\end{document}